\documentclass[10pt]{article}
\usepackage{subfiles}
\usepackage{geometry}
\usepackage{float}% pour le placement des figure
\usepackage{setspace}
\usepackage{url}
\usepackage{enumitem}
\usepackage{amstext}
\usepackage{xcolor}
\usepackage{amsfonts}
\usepackage{amssymb}
\usepackage{graphics}
\usepackage{graphicx}
\usepackage{pstricks,pst-node,pst-tree}
\usepackage{enumerate}
\usepackage{tikz} 
\usepackage{cutwin}
\usepackage{xcolor}
\usepackage{colortbl}
\usepackage{algorithm,algorithmic}
\usepackage{multirow}
\usepackage{amsthm}
\usepackage{glossaries}
\usepackage{array,arydshln}
\usetikzlibrary{shapes,matrix}
\usepackage{comment}
\usepackage{mathrsfs}
\usepackage{authblk}
\usepackage{wrapfig}
\usepackage[font=footnotesize]{caption}
\usepackage{stackengine}
\usepackage{caption}
\usepackage{subcaption}
\usepackage{hyperref}
\usepackage{stmaryrd}
\usepackage{bbm}
\usepackage{relsize}

\usepackage{xr}
\makeatletter

\newcommand*{\addFileDependency}[1]{% argument=file name and extension
\typeout{(#1)}% latexmk will find this if $recorder=0
\@addtofilelist{#1}
\IfFileExists{#1}{}{\typeout{No file #1.}}
}\makeatother

\newtheorem{theorem}{Theorem}
\newtheorem{proposition}{Proposition}
\newtheorem{corollary}{Corollary}

\newtheorem{remark}{Remark}

\newcommand  {\RR}{\mathcal{R}}
\newcommand  {\GG}{\mathcal{G}}

\newcommand{\noni}[2]{\overline{#1}^{#2}}

\newcommand{\N}{\mathbb{N}}

\definecolor{sable}{RGB}{153,51,0}
\definecolor{macouleur}{RGB}{180,0,0}
\definecolor{ForestGreen}{RGB}{0,180,0}
\definecolor{Aquamarine}{rgb}{0.5, 1.0, 0.83}
\definecolor{aqua}{rgb}{0.0, 1.0, 1.0}
\definecolor{darkcyan}{rgb}{0.0, 0.55, 0.55}
\definecolor{darkpastelpurple}{rgb}{0.59, 0.44, 0.84}
\definecolor{ao(english)}{rgb}{0.0, 0.5, 0.0}
\definecolor{atomictangerine}{rgb}{1.0, 0.6, 0.4}
\definecolor{amaranth}{rgb}{0.9, 0.17, 0.31}

\title{On a theorem of François Robert}

\author[1]{Brigitte Mossé}
\author[2]{Sasha Pignol}
\author[1]{Elisabeth Remy}
\affil[1]{Aix Marseille Univ, CNRS, I2M (UMR 7373) Marseille, France}
\affil[2]{Aix Marseille Univ, Marseille, France}
\begin{document}
\maketitle

\vspace{1cm}
\paragraph{Abstract} A well-known theorem by François Robert expresses the degenerated character of a synchronous Boolean finite dynamical system, in the case where the associated regulatory graph does not contain any circuit: all states of the system go towards a single fixed point. We present a large family of updating modes of Boolean models with the same particularity.
\newpage
\tableofcontents

%%%%%%%%%%%%%%%%%%%%%%%%%%%  
%%%%%%%%%%%%%%%%%%%%%%%%%%%
%%%%%%%%%%%%%%%%%%%%%%%%%%%
\section{Introduction}

In the field of finite dynamical systems, links between interaction graphs and dynamics are widely studied. They are particularly useful for certain applications, in particular in biology. For example, models of genetic interaction networks generate spaces of states and trajectories whose size leads to a combinatorial explosion, preventing the analysis of the dynamical properties of the system.  The study of topological properties of genetic regulation networks is more accessible than the analysis of the dynamics, and provides access to certain characteristics of the corresponding dynamics.
A theorem of François Robert states that, under the hypothesis that a regulatory graph does not contain any circuit, the synchronous dynamics of a related Boolean model is reduced to trajectories towards a single stable state \cite{robert86a, robert95}.
Hence, the “complexity” of the dynamics comes from the presence of circuits in the regulatory graph. In biology, the regulatory graphs are directed signed graphs, the sign representing the nature of the regulation (activation or inhibition).
Regulatory circuits, also called feedback loops, are already well-known to play significant dynamical roles. Positive regulatory circuits (with an even number of negative interactions) have been associated with multi-stability, which may account for biological differentiation phenomena. 
Negative regulatory circuits (with an odd number of negative interactions) have been associated with sustained periodic behaviors as homeostasis \cite{thomas81a}.
Mathematical results proving such properties, so called Thomas’s rules, emphasize their links with the complex notion of circuit functionality \cite{remy08a,richard10a}.
Robert's theorem is a general result linking topological property of the regulatory graph and characteristics of the associated dynamics. In \cite{robert86a, robert95}, F. Robert extends the theorem to other updating modes than the usual synchronous one (the Gauss-Seidel type of synchronous updating) and suggests that there are many generalizations. We give a series of situations for which this is indeed true.

%%%%%%%%%%%%%%%%%%%%%%%%%%%  
%%%%%%%%%%%%%%%%%%%%%%%%%%%
%%%%%%%%%%%%%%%%%%%%%%%%%%%

\section{Basics on Boolean finite dynamical systems }

A {\it Boolean model} is a map $S: X=\{0,1\}^n \to X$, with $n$ an integer $>0$. The elements of $X$ are called {\it states}. 
For $x \in X$ and a non-empty part $J$ of $\{0,1\}^n$, we will denote by $\noni{x}{J}$ the state obtained from $x$ by switching coordinates with index in $J$, and $\noni{x}{i}=\noni{x}{\{i\}}$ if $J$ is a singleton.

A {\it regulatory graph}, denoted by $\mathcal{R}\mathcal{G}(S)$, is classically associated to a Boolean model $S$. The $n$ vertices $g_1$, ... , $g_n$ of $\RR \GG(S)$ are abstract components whose {\it level} is given by the $x_i$ along the dynamics. There is an edge in $\mathcal{R}\mathcal{G}(S)$ from $g_i$ to $g_j$, called {\it interaction} from $g_i$ to $g_j$, if there exist two states $x$ and $y=\noni{x}{i}$ such that $S_j(x) \neq S_j(y)$. This condition means that $g_i$ is a regulator of $g_j$: the level $x_i$ of $g_i$ has an influence on the level $x_j$ of $g_j$ under the action of $S$, influence expressed on at least one pair of states $\{ x, \noni{x}{i} \}$ and interpreted as an activation or an inhibition, as the case may be.

Given the map $S=(S_1, \dots,S_n)$, to any state $x$ is associated its {\it updating set}
$Upd_S(x)=\{ i \in \{1, \dots n\}\,; S_i(x) \neq x_i \}$. This brings to consider the most usual finite dynamical system associated to $S$, that is the {\it synchronous dynamics}, with a simultaneous change of all the coordinates of the state $x$ specified by $Upd_{S}(x)$. Any state has exactly one successor, its image under $S$, so that this dynamics is deterministic. The graph of the map $S$ is called synchronous {\it state transition graph} (STG). This graph expresses the iteration of $S$, and its paths give the trajectories of the dynamics.

Two other updating modes are of interest, leading to the {\it asynchronous} and {\it fully asynchronous} dynamics \cite{garg08a,pauleve:hal-03736899}. Given a state $x$, the asynchronous STG has edges of origin $x$ those of extremity $\noni{x}{\{i\}}$, where $i$ is any element of $Upd_S(x)$.
The fully asynchronous STG has edges of origin $x$ those of extremity $\noni{x}{J}$ where $J$ is any non-empty part of $Upd_S(x)$. In both cases, the trajectories are all the paths of the STG. Both dynamics are non-deterministic, with states potentially having several successors. Their interest in applications is to simultaneously take into account various possibilities of delays or priorities to be given to switching.

%%%%%%%%%%%%%%%%%%%%%%%%%%%  
%%%%%%%%%%%%%%%%%%%%%%%%%%%
%%%%%%%%%%%%%%%%%%%%%%%%%%%
\section{Robert's theorem}

We recall Robert's theorem, and main ingredients of the proof given in \cite{robert86a, robert95}.

\medskip
Let us consider a Boolean model $S$, any given updating mode, and the dynamics and STG they determine. These dynamics and STG are said \textit{simple} if the STG has only one attractor (i.e. one terminal strongly connected component), which is reduced to a single state.

\begin{theorem}[Robert's theorem] \label{THE theorem}
    Let $S$ be a Boolean model on $X = \{0,1\}^n$. If $\mathcal{RG}(S)$ does not contain any circuit, then the synchronous dynamics defined by $S$ is simple. Moreover, for all $x\in X$, the sequence $(S^k(x))_{k\in\N}$ converges in at most $n$ iterations to the single fixed point of $S$.
\end{theorem}

\begin{proof}

In \cite{robert86a, robert95}, Robert adopts a metric and algebraic point of view to show this result; we give a sketch of this proof. 
The framework is that of Boolean calculation (that is $0.0=0.1=1.0=0$, $1.1=1$, $0+0=0$ and $0+1=1+0=1+1=1$). A \textit{vector Boolean distance} between states is defined coordinate by coordinate by $d(x,y)=(\delta(x_1,y_1),\dots,\delta(x_n,y_n))$, where 
$\delta(a,b)= \lvert a-b \lvert$. Inequalities are also defined coordinate by coordinate. The proof is then based on two arguments.

The first one is a basic inequality, consequence of the triangular inequality for $d$. A matrix $B(S)$ is defined as the $n\times n$ Boolean matrix transpose of the adjacency matrix of $\mathcal{RG}(S)$, i.e.  $b_{ij} = 1$ if there is an edge in $\mathcal{RG}(S)$ from $g_j$ to $g_i$, and $b_{ij} = 0$ otherwise. Then, for all $x,\;y \; \in X$ one has 
    $d(S(x),S(y))^t \; \leq \; B(S)\;d(x,y)^t$, where "$. ^t$" denotes the matrix transposition.

The second argument is the translation on $B(S)$ of the assumption that $\mathcal{RG}(S)$ does not contain any circuit. 
Let us recall that if $G = (V,E)$ is a directed graph with $N$ vertices such that $G$ does not contain any circuit, a topological sorting of $G$ is any process of assigning to each vertex $v$ a number $i\in \{1, \dots N\}$, such that for any $j\in \{1, \dots N\}$, the inequality $i\le j$ implies $(v_j,v_i)\notin E$ \cite{Bang-Jensen}. The existence of such a process implies conversely that $G$ is circuit-free. The hypothesis on $\mathcal{RG}(S)$ is equivalent to the existence of a permutation matrix $P$ such that $P^tB(S)P $ is strictly lower triangular. Indeed, such a conjugation of the matrix $B(S)$ amounts to a renumbering of the components by $P$ performing a topological sort. 

The second argument implies that $B(S)^n=0_{\mathcal{M}_n(\{0,1\})}$, and iterating the basic inequality gives the conclusion. \end{proof}

In \cite{robert86a, robert95}, Robert shows that the updating modes of type Gauss-Seidel lead to new synchronous versions of this theorem. Let us give the definition of the simplest of these modes.
 Given the Boolean model $S=(S_1, \dots,S_n)$ and $x=(x_1, \dots,x_n) \in X$, the Gauss-Seidel Boolean model $G=(G_1, \dots,G_n)$ associated to $S$ is obtained by the following iteration: 
 $G_1(x)=S_1(x)=y_1$, $G_2(x)=S_2(y_1, x_2, \dots,x_n)=y_2$, ... , $G_n(x)=S_n(y_1, y_2, \dots,y_{n-1}, x_n)$.
  Then if $S$ satisfies the hypothesis of Theorem \ref{THE theorem}, the synchronous STG of $G$ is simple, and states reach the single fixed point of $S$ in at most $n$ iterations of $G$.

  \bigskip
  We now give other versions of the theorem. The following fully asynchronous one is an immediate consequence of Theorem \ref{THE theorem}.

\begin{corollary}
    Let $S$ be a Boolean model on $X = \{0,1\}^n$. If $\mathcal{RG}(S)$ does not contain any circuit, then the fully asynchronous dynamics defined by $S$ is simple. Moreover, for all $x\in X$, there exists in the fully asynchronous STG a path of length at most $n$ from $x$ towards the single fixed point of $S$.
\end{corollary}

   \begin{remark} 
In the same spirit as the fully asynchronous version, we can look at the case of the Most Permissive updating mode, which introduces intermediate levels between 0 and 1 and leads to a significant increase in reachability between Boolean states. We refer the reader to \cite{pauleve_reconciling_2020} for a precise definition. It is immediate to validate a Most Permissive version of Robert's theorem. Indeed, in the Most Permissive STG of $S$, every non-Boolean state goes to a Boolean state, and from every Boolean state starts its synchronous trajectory, inserting if necessary non-Boolean states. So under hypothesis of Theorem \ref{THE theorem}, any state reaches the single fixed point of the synchronous dynamics in a path involving at most $n+1$ Boolean states.
   \end{remark} 
  
   \begin{remark} 
   Under the hypothesis of Theorem \ref{THE theorem}, the regulatory graph has at least one component without regulator. If it is the case of $g_i$, that means that the function $S_i$ is constant (see details in Proposition \ref{prop:constant}). One might be more interested in models with instantiated inputs, that are components $g_i$ for which $S_i(x)=x_i$ for all $x\in X$. Remark that the inputs are self-activated and have no regulator but themselves. To give an instantiation to input $g_i$ consists in considering the part of the space such that its level is equal to a fixed element $a_i$ of $\{0,1\}$. It is easy to give a variant of Robert's theorem in this situation, which we briefly present.
\end{remark}

\begin{theorem}
    Let $S$ be a Boolean model on $X = \{0,1\}^n$. We suppose that the model admits $r$ inputs $g_1, \dots , g_r$, with $1\leq r\leq n$. For any $a=(a_1, \dots ,a_r) \in \{0,1\}^r$, let us denote $\mathcal{C}_a=\{x\in X \,;\, \forall i \in \{1, \dots,r\},\; x_i=a_i\}$. If $\mathcal{RG}(S)$ does not contain any circuit except the self-activation of the inputs, then the synchronous STG of $S$ admits as only attractors $2^r$ fixed points, determined by their $r$ first coordinates. For each input instantiation $a$, the set $\mathcal{C}_a$ is the basin of attraction of a single fixed point, reached from any element of $\mathcal{C}_a$ in at most $n-r$ iterations of $S$.
\end{theorem}

We point out that for a detailed proof, it is enough to adapt the basic inequality and matrix considerations of the proof of Theorem \ref{THE theorem} by restrictions to sets $\mathcal{C}_a$.

%%%%%%%%%%%%%%%%%%%%%%%%%%%  
%%%%%%%%%%%%%%%%%%%%%%%%%%%
%%%%%%%%%%%%%%%%%%%%%%%%%%%
\section{Asynchronous version of Robert's theorem}

\begin{figure}[h]
  \begin{subfigure}[b]{.25\linewidth}
    \centering
    \begin{tabular}{ c|c }
      $x$ & $S(x)$ \\
      \hline
      00 & 11 \\
      01 & 00 \\
      10 & 00 \\
      11 & 11
    \end{tabular}

    \caption{Boolean model $S$}
  \end{subfigure}\hfill
  \begin{subfigure}[b]{.25\linewidth}
    \centering
    \includegraphics[width=0.7\linewidth]{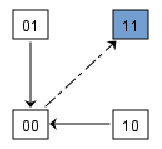}
    \caption{Synchronous STG of $S$ - the graph is simple.}
  \end{subfigure}\hfill
  \begin{subfigure}[b]{.25\linewidth}
    \centering
    \includegraphics[width=0.7\linewidth]{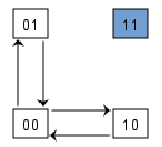}
    \caption{Asynchronous STG of $S$ - the graph is not simple.}
  \end{subfigure}

  \caption{Example of Boolean model with simple synchronous STG and not simple asynchronous STG}
\label{fig:counterexample}
\end{figure}

The first remark, considering asynchronous dynamics, is that an asynchronous version of Robert's theorem is not an immediate corollary of Theorem \ref{THE theorem}.
Indeed, for a given Boolean model $S$, the fact that the synchronous STG of $S$ is simple does not guarantee that its asynchronous STG will also be simple.
The STGs depicted on figure \ref{fig:counterexample} provide a counterexample to this, for $S: \{0,1\}^2 \to \{0,1\}^2$ given by $ S_1(x)=S_2(x)=(x_1\land x_2)\lor(\lnot x_1\land\lnot x_2)$.

In order to develop a proof, we will consider subsets of $X$ of states whose some coordinates are fixed. Given $I$ a non-empty subset of $\{1,\dots,n\}$ and $a=(a_i)_{i\in I}$ a family of elements of $\{0,1\}$, we set $\mathcal{C}_a = \{x\in X\;;\; \forall i \in I,\; x_i=a_i\}$.

\begin{proposition}\label{prop:constant}   Let $S$ be a Boolean model on $X = \{0,1\}^n$. \begin{itemize}
    \item[--] For any $i\in \{1,\dots,n\}$, if $S_i$ is not constant, then the component $g_i$ is regulated by at least a component $g_j$ with $j\in \{1,\dots,n\}$ (we may have $i=j$). 
\item[--] Let $I$ be a non-empty subset of $\{1,\dots,n\}$ and $a=(a_i)_{i\in I}$ a family of elements of $\{0,1\}$.
For any $i\in I$, if $S_i$ is not constant on $\mathcal{C}_a$, then the component $g_i$ is regulated by at least a component $g_j$ with $j\in I$ (we may have $i=j$). 
\end{itemize} 
\end{proposition}

\begin{proof} \begin{itemize}
 \item[--] Let $i\in \{1,\dots,n\}$, and suppose that $S_i$ is not constant on $X$. Let us show that there exists $x\in X$ such that $S_i(x)\neq S_i(\overline x^j)$, which means that $g_i$ is regulated by $g_j$. Indeed, there exist $y,z\in X$ such that $S_i(y)\neq S_i(z)$. Then, one of the following non-equalities is true:\\
         $S_i(y_1,~ \dots,~ y_n) \neq S_i(z_1,~y_2,~ \dots,~ y_n)$,\\
        $S_i(z_1,~y_2,~ \dots,~ y_n) \neq S_i(z_1,~z_2,~y_3,~ \dots,~ y_n)$,\\
        $~~ \vdots$\\
        $S_i(z_1,~ \dots,~z_{n-1} ,~ y_n) \neq S_i(z_1,~ \dots,~ z_n)$,\\
      from which follows the existence of a suitable state $x\in X$.
    \item[--] Let $I$ be a non-empty subset of $\{1,\dots,n\}$ and $a=(a_i)_{i\in I}$ a family of elements of $\{0,1\}$. Suppose that $S_i$ is not constant on $\mathcal{C}_a$ for some $i\in I$. Then we proceed as above, with $y,z\in \mathcal{C}_a$ such that $S_i(y)\neq S_i(z)$. All the states appearing in the sequence of non equalities that we consider belong to $\mathcal{C}_a$, that gives the conclusion.
         
     \end{itemize}
\end{proof}

\begin{theorem} \label{async theorem}
 Let $S$ be a Boolean model on $X = \{0,1\}^n$. If $\mathcal{RG}(S)$ does not contain any circuit, then the asynchronous dynamics defined by $S$ is simple.
    Moreover, for any $x\in X$, there exists in the asynchronous STG a path of length at most $n$ from $x$ towards the single fixed point of $S$.
\end{theorem}

\begin{proof}

 Assume that $\mathcal{RG}(S)$ does not contain any circuit, and let us rename components, using a topological sort.
 
   - The component $g_1$ is not regulated by any component.
    Thus, by Proposition \ref{prop:constant}, the map $S_1$ is constant on $\mathcal{C}_0=X$. Let us denote by $\alpha_1$ its single value.
    
    - Let $\mathcal{C}_1=\{x\in X\;;\;  x_1=\alpha_1\}$. The component $g_2$ may only be regulated by the component $g_1$. Thus, by Proposition \ref{prop:constant}, $S_2$ is constant on $\mathcal{C}_1$. Let us denote $\alpha_2$ its single value on $\mathcal{C}_1$.
    
    Iterating this process, we get that for $k\in \{2,\dots,n\}$, the map $S_k$ is constant on the set $\mathcal{C}_{k-1}=\{x\in X\;;\; \forall i \in \{1, \dots,k-1\},\; x_i=\alpha_i\}$, with single value $\alpha_k$ on $\mathcal{C}_{k-1}$. The sets $\mathcal{C}_k$ are of cardinality $2^{n-k}$, and finally the image of $\mathcal{C}_{n-1}$ by the map $S_n$ is reduced to a single state $\alpha=(\alpha_1, \dots,\alpha_n)$.
    
    \medskip
By construction, the state $\alpha=(\alpha_1, \dots,\alpha_n)$ is a fixed point of $S$, thus it is the single fixed point of $S$ given by Theorem \ref{THE theorem}. Moreover, if one starts from an arbitrary state $x\in X$, successively applying $S_1$, $S_2$, ... ensures to obtain an asynchronous path from $x$ towards $\alpha$ passing through $\mathcal{C}_1$, $\mathcal{C}_2$, ... of length at most $n$.

\end{proof}

%%%%%%%%%%%%%%%%%%%%%%%%%%%  
%%%%%%%%%%%%%%%%%%%%%%%%%%%
%%%%%%%%%%%%%%%%%%%%%%%%%%%
\section{Extension to a family of updating modes}

Starting from a Boolean model $S$ on $X=\{1,\dots,n\}$, we now introduce a family of Boolean dynamics encompassing those we have already mentioned.
We consider a set $\mathcal{P}$ of non-empty parts of $\{1,\dots,n\}$ of union $\{1,\dots,n\}$.
        Let us associate with $\mathcal{P}$ the updating mode whose STG, denoted STG$_{\mathcal{P}}$, has for edges with origin a state $x$ those of extremities $\noni{x}{J \cap Upd_S(x) }$, where $J$ is any element of $\mathcal{P}$ such that $J \cap Upd_S(x) \neq \emptyset$. 
         Thus, the synchronous mode corresponds to the case where the only element of $\mathcal{P}$ is $\{1,\dots,n\}$, the asynchronous mode to the case where the singletons $\{i\}$, with $1\leq i\leq n$, are the elements of $\mathcal{P}$, and the fully asynchronous mode to the case where the $2^n-1$ non-empty parts of $\{1,\dots,n\}$ are the elements of $\mathcal{P}$. The assumption on $\mathcal{P}$, whose elements cover $\{1,\dots,n\}$, ensures that no information contained in $S$ is lost.
         We will call this type of dynamics \textit{$\mathcal{P}$-asynchronous dynamics}.

         Such an extension may be useful for modeling interaction networks subject to certain priorities, or collaborations between components \cite{faure06a}.

  \begin{theorem} \label{theorem 4}
 Let $S$ be a Boolean model on $X = \{0,1\}^n$ and $\mathcal{P}$ a set of non-empty parts of $\{1,\dots,n\}$ of union $\{1,\dots,n\}$. If $\mathcal{RG}(S)$ does not contain any circuit, then the $\mathcal{P}$-asynchronous dynamics defined by $S$ is simple.
    Moreover, for any $x\in X$, there exists in STG$_{\mathcal{P}}$ a path of length at most $n$ from $x$ towards the single fixed point of $S$.
\end{theorem}

\begin{proof} The proof is just an adaptation of the one of Theorem \ref{async theorem}. We first repeat this proof word for word until we obtain the fixed point $\alpha=(\alpha_1,\dots,\alpha_n)$. We then select a sequence $P_1, P_2, \dots , P_n$ of elements of $\mathcal{P}$, possibly including repetitions, such that $1\in P_1$, $2\in P_2$, ..., $n\in P_n$. Starting now from an arbitrary state $x\in X$, successively applying synchronously $\{S_i\;;\; i\in P_1\}$, $\{S_i\;;\; i\in P_2\}$, ... ensures to obtain an $\mathcal{P}$-asynchronous path from $x$ towards $\alpha$ passing through $\mathcal{C}_1$, $\mathcal{C}_2$, ... of length at most $n$.

\end{proof} 
     \begin{remark}
This proof thus covers proofs of the synchronous and asynchronous cases.
     \end{remark}

\begin{remark}
    We could imagine other types of updating modes, leading to a pruning of state transition graphs, and guided for example by an application framework. The following proposition shows that under the hypothesis of Robert's theorem, it is impossible to create other attractors than stable points by pruning.
\end{remark}

            \begin{proposition} Let $S$ be a Boolean model on $X = \{0,1\}^n$ and $\mathcal{P}$ a set of non-empty parts of $\{1,\dots,n\}$ of union $\{1,\dots,n\}$. If $\mathcal{RG}(S)$ does not contain any circuit, then STG$_{\mathcal{P}}$ does not contain any cycle of length $\geq 2$.
         \end{proposition}

\begin{proof}  Assume that $\mathcal{RG}(S)$ does not contain any circuit, and let us rename components, using a topological sort.  
Suppose that there exist two distinct states $x$ and $y$ contained in a cycle $\gamma$ of STG$_{\mathcal{P}}$. Let $r=\inf \{i\in \{1, \dots,n\}\;;\; x_i\neq y_i\}$.
Proposition \ref{prop:constant} then gives the following. 

- If $r=1$, then at least one element of $\mathcal{P}$ involved in a transition of $\gamma$ contains $1$. But once $S_1$ acts on a state, all consecutive states have the same first coordinate. This is in contradiction with the fact that $x$ is consecutive to $y$, and vice-versa, thus $r\geq 2$.

- If $r=2$, then at least one element of $\mathcal{P}$ involved in a transition of $\gamma$ contains $2$. The component $g_2$ may be regulated or not by $g_1$, but the map $S_2$ is constant on each of the two subspaces $\{x\in X\;;\; x_1=0\}$ and $\{x\in X\;;\; x_1=1\}$. As above, this implies that all the states of $\gamma$ have the same second coordinate, in contradiction with $r=2$, thus $r\geq 3$.

- If $r=3$, then at least one element of $\mathcal{P}$ involved in a transition of $\gamma$ contains $3$. The component $g_3$ may be regulated or not by $g_1$ and by $g_2$, but the map $S_3$ is constant on each of the four subspaces $\{x\in X\;;\; (x_1,x_2)=(a_1,a_2)\}$, where $a_1$ and $a_2$ are equal to $0$ or $1$. Thus all the states of $\gamma$ have the same third coordinate, in contradiction with $r=3$, and so on, hence the result.

\end{proof}

\section{Conclusion}

This work establishes a generalization of Robert's theorem, encompassing both synchronous and fully synchronous updating modes of Boolean models. This advancement opens new research perspectives, particularly regarding the possibility of extending these results whatever the Boolean network updating modes. In this context, a natural and promising extension would be the exploration of the multi-valued framework. This approach holds particular importance in the modeling of biological graphs, where system complexity often requires a more nuanced representation than the simple binary framework. A differential version of the theorem has already given rise to a development, in the study of the role of circuits in dynamical systems described by systems of differential equations developped in \cite{kaufman07a}. Robert's theorem seems to be an intrinsic result for models whose regulation graphs are circuit-free, regardless of the framework, and this question could be investigated.

\vspace{1cm} All these developments and results were initially the subject of a Bachelor internship report by Sasha Pignol, as part of his mathematical course at Aix Marseille University, supervised by the other authors of the article.

\bibliographystyle{abbrv}

\end{document}